\documentclass[12pt,oneside,a4paper]{article}
\usepackage[english]{babel}
\usepackage[utf8]{inputenc}
\usepackage[T1]{fontenc}
\usepackage{latexsym}
\usepackage{graphicx}
\usepackage{amsmath}
\usepackage{amsfonts}
\usepackage{amssymb}
\usepackage{amsthm}
\author{Julien Cassaigne and Idrissa Kabor\'e}
\title{Word of low complexity without uniform frequencies}
\date{\ }
\newtheorem{definition}{Definition}[section]
\newtheorem{thm}{Theorem}[section]
\newtheorem{prop}[thm]{Proposition}

\newtheorem{lm}{Lemma}[section]

\newtheorem{notation}{Notation}[section]

\newcommand{\N}{\mathbb{N}}

\newcommand{\ie}{\emph{i.e.}}

\begin{document}
\maketitle

\textbf{Abstract:} In this paper, we construct a uniformely recurrent infinite word of low complexity without uniform frequencies of letters. This shows the optimality of a bound of Boshernitzan, which gives a  sufficient condition for a uniformly recurrent infinite word to admit uniform frequencies.

\bigskip
 
\emph{Mathematics Subject Clasification: 37B10, 37A25, 68R15 }.

\section{Introduction}
Let us consider an infinite word $u$ over a finite alphabet. We can naturally associate to it a subshift. The goal of this paper is to describe some ergodic properties of this subshift. By Oxtoby theorem, we know that the subshift is uniquely ergodic if and only if, in $u$, each finite word has uniform frequency. Moreover the subshift is minimal if $u$ is uniformly recurrent. 

For a long time, people have tried to find some  conditions on infinite words which imply one of these properties. M. Keane gave in \cite{keane} a uniformly recurrent infinite  word with complexity $3n+1$  (from 4-interval exchange map) which does not possess uniform frequencies. Later, Boshernitzan, in \cite{bosher}, proved that a uniformly recurrent infinite word admits uniform frequencies if either of the following sufficient conditions is satisfied: \[\liminf\dfrac{\mathbf{p}(n)}{n}<2\ \textrm{or}\ \limsup\dfrac{\mathbf{p}(n)}{n}<3\]
where $\mathbf{p}$ denotes the complexity function of $u$ (see \cite{cant} chap. 4, by J. Cassaigne and F. Nicolas  for more details on the subject).

The bound of the second sufficient condition of Boshernitzan being already optimal by  Keane's result, our goal in his work is to establish the optimality of the bound of the first sufficient condition. And we succeed to  construct a  uniformly recurrent infinite word without uniform frequencies,  the complexity function of which verifies $\liminf\dfrac{\mathbf{p}(n)}{n}=2$.

This result relates some properties of the complexity function and the ergodic measures of the subshift. This type of question has been  investigated in the last years. The goal is to bound the number of ergodic measures of the subshift in terms of the complexity function.

Boshernitzan was the first to look at it, see \cite{bosher}. During his Phd. thesis T. Monteil, see \cite{monteil} and \cite{cant} (chap. 3 by S. Ferenczi and Th. Monteil) has proved, the same result with different techniques: If $\limsup\dfrac{\mathbf{p}(n)}{n}=K\geq 2$, then the subshift has at most $K-2$ ergodic measures. 
Since this time, some results have appeared in the same veine: V. Cyr and B.  Kra have also obtained similar results, see \cite{cyrkra1, cyrkra2}. In the first paper, they prove that the bound of Boshernitzan is sharp. In the second one, they construct minimal subshifts with complexity function arbitrarily close to linear but having uncountably many ergodic measures. We can also cite M. Damron and D. Fickenscher \cite{dafi} who obtained the bound $\frac{K+1}{2}$ under a condition on the bispecial words.

Nevertheless, it seems that our proof is of a different nature, with an explicit construction of the infinite word.

After the preliminaries (section 2) we construct an infinite word which is uniformly recurrent in  section 3, then we show in section 4 that this word is without uniform frequencies of letters and to finish we study the complexity of this word in section 5 and give in section 6 the proof of the main statement of section 3.

\section{Preliminaries}

In all that follows we consider  the alphabet $\mathcal{A}=\left\lbrace 0,\,1 \right\rbrace $.  Let us denote $\mathcal{A}^\ast$,   the set of  the finite words on alphabet $\mathcal{A}$, $\varepsilon$ the empty word. For all $u$ in $\mathcal{A}^\ast$, $|u|$ denotes the length (the number of letters it contains) of the word $u$ ($\left|\varepsilon\right|=0$) and for any letter $x$ of $\mathcal{A}$, $|u|_x$ is the number of occurrences of the letter $x$ in $u$. We call Parikh vector  of a finite word $u$,    the vector denoted by $U$ and defined by $\begin{pmatrix}\left|u\right|_0\\\left|u\right|_1
\end{pmatrix}.
$

A finite word $u$ of length $n$ formed by repeating a single letter $x$ is typically denoted $x^n$.  We define the $n$-th power of a finite word $w$ as being the concatenation  of $n$ copies of $w$; we denote it $w^n$. An infinite word is an infinite sequence of letters of $\mathcal{A}$. We denote $\mathcal{A}^\omega$ the set of infinite words on $\mathcal{A}$.

We say a  finite word $v$ is a factor of  $u$  if there exist two words $u_1$ and $u_2$  on the alphabet $\mathcal{A}$ such that $u=u_1vu_2$; we also say  that $u$ contains $v$. The factor $v$ is said prefix (resp. suffix) if $u_1$ (resp. $u_2$) is the empty word. For any word $u$, the set of  factors  of length $n$ is denoted $\mathcal{L}_n(u)$. The set of all factors of $u$ is simply denoted $\mathcal{L}(u)$. 

\begin{definition}
 Let $u$ be an infinite word on the alphabet $\mathcal{A}=\left\lbrace0,\ 1\right\rbrace  $. A factor $v$ of $u$ is said to be
\begin{itemize}

\item  a right special  factor if $v0$ and $v1$ are both factors of $u$, and a left special factor if $0v$ and $1v$ are both factors of $u$.

\item  a bispecial factor of $u$ if $v$ is simultaneously a right special factor and a left special factor of $u$.

\item  a stronq bispecial factor of $u$ if $0v0,\ 0v1,\ 1v0,\ 1v1$ are factors of $u$ and a weak bispecial factor if uniquely $0v0$ and  $1v1$, or $0v1$ and $1v0$, are factors of $u$.

\item an ordinary bispecial factor of $u$ if $v$ is a bispecial factor of $u$ which is neither strong nor weak. \end{itemize}

\end{definition}

An infinite word $u$ is said to be recurrent if any factor of $u$ appears infinitely often. An infinite word $u$ is uniformly recurrent   if for all $n\in \N$, there exists $N$ such that any factor of $u$ of length $N$ contains all the factors of $u$ of length $n$.

\begin{definition}
 Let $u$ be an infinite word on an alphabet $\mathcal{A}$. The complexity function of $u$ is a function counting the number of distinct factor of $u$ of length $n$  for any given  $n$. It is denoted $\mathbf{p}$ and so that:
\[ \mathbf{p}(n)=\#\mathcal{L}_n(u).\]
\end{definition}

Let us denote $\mathbf{s}$ and $\mathbf{b}$ the functions respectively called first difference and second difference of the complexity of $u$; they are defined  as follows: $\mathbf{s}(n)= \mathbf{p}(n+1)-\mathbf{p}(n)$ and  $\mathbf{b}(n)=\mathbf{s}(n+1)-\mathbf{s}(n)$.

On a binary alphabet the function $\mathbf{s}$ counts the number of special factors for a given length in $u$. Let us denote $\mathbf{m}$ the map from $\mathcal{L}(u)$ to $\left\lbrace -1,\, 0,\,+1\right\rbrace $ defined by
$$\forall v\in \mathcal{L}(u),\ \mathbf{m}(v)=\left\{
\begin{array}{ll}
        -1& \textrm{if}\  v\ \textrm{is}\ \textrm{weak}\ \textrm{bispecial}\\
        +1&  \textrm{if} \  v\ \textrm{is}\ \textrm{strong}\ \textrm{bispecial}\\
        0&  \textrm{otherwise}\\
	
\end{array}\right.$$
The   following formula  was given by the first author in \cite{cassaigne2}:
\[\forall n\geq 0, \ \mathbf{s}(n)=1+\sum _{\begin{array}{c}w\ \in \mathcal{L}(u)\\
                                            |w|<n
                                           \end{array}}\mathbf{m}\left( w\right)=1+\sum _{\begin{array}{c}w\ \textrm{bispecial}\\
                                            |w|<n
                                           \end{array}}\mathbf{m}\left( w\right) .
\]

This relation allows  to compute the complexity $\mathbf{p}(n)$  provided  when we are able to describe the set of strong and weak bispecial factors of the binary infinite word $u$.
 
\begin{definition}
 Two bispecial factors $v$ and $w$  of an infinite word $u$ on the alphabet $\left\lbrace 0,\ 1\right\rbrace $ are said to have the same type if they are all strong, weak, or ordinary. In other words the bispecial $v$ and $w$ have the same if $\mathbf{m}(v)=\mathbf{m}(w)$. 
\end{definition}

\begin{definition} (\cite{cant} chap. 7, by S. Ferinczi and T. Monteil) Let $u$ be an infinite word on an alphabet $\mathcal{A}$.
 \begin{itemize}
  \item We say that  $u$ admits frequencies if for any factor $w$ , and any sequence $(u_n)$ of prefixes of $u$ such that $\lim_{n\rightarrow \infty}=\infty$, then $\lim_{n\rightarrow\infty}\frac{|u_n|_w}{|u_n|}$  exists. 
\item We say that  $u$ admits uniform frequencies if for any factor $w$, and any sequence $(u_n)$ of factors of $u$ such that $\lim_{n\rightarrow \infty}=\infty$, then $\lim_{n\rightarrow\infty}\frac{|u_n|_w}{|u_n|}$  exists. 
 \end{itemize}
\end{definition}

In \cite{keane}, M. Keane gave  an example of a uniformly recurrent infinite word with complexity $3n+1$ which does not possess uniform frequencies. Later, Boshernitzan \cite{bosher} obtained the following results:

\begin{thm} Let $u$ be an infinite word on an alphabet $\mathcal{A}$. Then, $u$  admits uniform frequencies if its complexity function verifies at least one of the  following conditions:
 \begin{itemize}
  \item   $\displaystyle{\liminf\dfrac{\mathbf{p}(n)}{n}<2},$
  \item   $\displaystyle{\limsup\dfrac{\mathbf{p}(n)}{n}<3}.$
 \end{itemize}
\end{thm}

The example of Keane enssures  that constant 3 is optimal in the second condition, \ie, it cannot be replaced with a larger constant.

\section{Construction of a class of  uniformly recurrent words}

Let $(l_i)$, $(m_i)$, $(n_i)$ be three integer sequences which are strictly increasing and verify the following conditions:

\begin{itemize}

\item[$\bullet$] $l_i<m_i<n_i$,

\item[$\bullet$] $\dfrac{m_i}{l_i}$ increases exponentially to $+\infty$,

\item[$\bullet$] $\dfrac{n_i}{m_i}$ increases exponentially to $+\infty$.

\end{itemize}

Let us define in $\mathcal{A}^\ast$ two sequences $(u_i)$ and $(v_i)$ in the following way:
$u_0=0$, $v_0=1$ and for all $i\in\mathbb{N}$, $u_{i+1}=u_i^{m_i}v_{i}^{l_i}$ and $v_{i+1}=u_i^{m_i}v_{i}^{n_i}$. The sequence $(u_i)$ converges towards an infinite word $u$.

For $i\geq 1$, consider the substitution $\sigma_i$ defined by $\sigma_i(0)=0^{m_i}1^{l_i}$, $\sigma_i(1)=0^{m_i}1^{n_i}$. Then, we have 

$u_{i}=\sigma_0\sigma_1\sigma_2\ldots \sigma_{i-1}(0)$ and $v_{i}=\sigma_0\sigma_1\sigma_2\ldots \sigma_{i-1}(1)$.

\begin{thm} \label{rec.unif}
 Any infinite word $u$ so defined is uniformly recurrent.
\end{thm}

The proof  is given at the end of this paper.

\section{The word u is without uniform frequencies}

\begin{lm} \label{produit} For all  $i\geq1$ we have:
\begin{enumerate}
\item $\displaystyle{\dfrac{\vert u_i\vert_0}{\vert u_i\vert}\geq \left(1+\dfrac{l_0}{m_0} \right)^{-1} \Pi_{j=1}^{i-1}\left(1+\dfrac{l_jn_{j-1}}{m_jl_{j-1}} \right)^{-1}}$ 
\item $\dfrac{\vert v_i\vert_1}{\vert v_i\vert}\geq \Pi_{j=0}^{i-1}\left(1+\dfrac{m_j}{n_j} \right)^{-1}
. $
\end{enumerate}
\end{lm}
 \begin{proof} $\bullet$ Lower bound on $\dfrac{\vert u_{i+1}\vert_0}{\vert u_{i+1}\vert}$.

Firstly, we have for all $i\geq 0$, $\vert u_{i}\vert\leq\vert v_i\vert$ since $\vert u_0\vert=\vert v_0\vert=1$ and $u_i$ is a strict prefix of $v_i$ for $i\geq1$. 
Then  $\dfrac{\vert v_i\vert}{\vert u_i\vert}=\dfrac{m_{i-1}\vert u_{i-1}\vert+n_{i-1}\vert v_{i-1}\vert}{m_{i-1}\vert u_{i-1}\vert+l_{i-1}\vert v_{i-1}\vert}\leq\dfrac{n_{i-1}}{l_{i-1}}$ since  $l_{i-1}< n_{i-1}$ for $i\geq1$.

As \[\vert u_{i+1}\vert_0=m_i\vert u_i\vert_0+l_i\vert v_i\vert_0\geq m_i\vert u_i\vert_0\]  and
\[ \vert u_{i+1}\vert=m_i \vert u_i\vert+l_i\vert v_i\vert=\vert u_i\vert\left( m_i+l_i\dfrac{\vert v_i\vert}{\vert u_i\vert}\right)\]
we deduce the following inequalities:
 \[\vert u_{i+1}\vert\leq \vert u_i\vert\left( m_i+l_i\dfrac{n_{i-1}}{l_{i-1}}\right) \]
 and
\[ \dfrac{\vert u_{i+1}\vert_0}{\vert u_{i+1}\vert}\geq \left( 1+\dfrac{l_i}{m_i}\dfrac{n_{i-1}}{l_{i-1}}\right)^{-1} \cdot\dfrac{\vert u_{i}\vert_0}{\vert u_{i}\vert}.\]
Thus 
\[\dfrac{\vert u_i\vert_0}{\vert u_i\vert}\geq \dfrac{\vert u_{1}\vert_0}{\vert u_{1}\vert}\Pi_{j=1}^{i-1}\left(1+\dfrac{l_jn_{j-1}}{m_jl_{j-1}} \right)^{-1}
. \]
$\bullet$ Lower bound on $\dfrac{\vert v_{i+1}\vert_1}{\vert v_{i+1}\vert}$.

We have \[\vert v_{i+1}\vert_1=m_i\vert u_i\vert_1+n_i\vert v_i\vert_1\geq n_i\vert v_i\vert_1\ \textrm{ and}\ \vert v_{i+1}\vert=m_i\vert u_i\vert+n_i \vert v_i\vert\leq\vert v_i\vert\left( m_i+n_i\right) \] since $\vert u_i\vert\leq\vert v_i\vert$.
So
\[\dfrac{\vert v_{i+1}\vert_1}{\vert v_{i+1}\vert}\geq \dfrac{n_i}{m_i+n_i}.\dfrac{\vert v_{i}\vert_1}{\vert v_{i}\vert}.\]
Hence
\[\dfrac{\vert v_{i}\vert_1}{\vert v_{i}\vert}\geq \Pi_{j=0}^{i-1}\left( 1+\dfrac{m_i}{n_i}\right) ^{-1}
.\] \end{proof}

In the rest  of the paper we need to fix

\begin{equation}
 l_i=2^{2.2^i+4},  m_i=2^{8.2^i}, \text{ and } n_i=2^{10.2^i}, \text{ for } i\geq 0. \tag{$\ast$}
\end{equation}

Then the inequalities of Lemma \ref{produit} become:
\begin{enumerate}
\item  $\displaystyle{\dfrac{\vert u_i\vert_0}{\vert u_i\vert}\geq \Pi_{j=1}^{i}\dfrac{1}{1+2^{-2^{j}}}}$
\item $\displaystyle{\dfrac{\vert v_i\vert_1}{\vert v_i\vert}\geq \Pi_{j=1}^{i}\dfrac{1}{1+2^{-2^{j}}}}$
\end{enumerate}
So we get
\begin{lm} \label{varianteproduit}
 \[\forall i\geq1,\ \min\left( \dfrac{\vert u_i\vert_0}{\vert u_i\vert},\,\dfrac{\vert v_i\vert_1}{\vert v_i\vert} \right)\geq \Pi_{j=1}^{i}\dfrac{1}{1+2^{-2^{j}}}.\]
\end{lm}

Then we have the following lemma:
\begin{lm} \label{excedent}\[\forall i\in\mathbb{N},\ \dfrac{\vert u_i\vert_0}{\vert u_i\vert}+\dfrac{\vert v_i\vert_1}{\vert v_i\vert}\geq \dfrac{3}{2}.\]
\end{lm}

 \begin{proof}
 \begin{itemize} \item[$\bullet$] For $i=0$, the inequality is evident.
                 \item[$\bullet$] For $i\geq1$, write: 
 $\ P_i=\Pi_{j=1}^{i}\dfrac{1}{1+2^{-2^j}}$. The sequence $(P_i)$ is decreasing and satisfies the following induction formula: $P_{i+1}=\dfrac{1}{1+2^{-2^{i+1}}}P_i$. 

Let us show, by induction, that  $\dfrac{4}{3}P_i=\dfrac{1}{1-2^{-2^{i+1}}}$.\\
We have  $\dfrac{4}{3} P_0=\dfrac{1}{1-2^{-2}}.$

 Assuming that for some $i\geq 0$, $\dfrac{4}{3}P_i=\dfrac{1}{1-2^{-2^{i+1}}}$ it follows : 
\[ \dfrac{4}{3}P_{i+1}=\dfrac{4}{3}P_i\times\dfrac{1}{1+2^{-2^{i+1}}}=\dfrac{1}{1-2^{-2^{i+1}}}\times\dfrac{1}{1+2^{-2^{i+1}}}=\dfrac{1}{1-2^{-2^{i+2}}}.
\]
So
\[P_i=\dfrac{3}{4}\times\dfrac{1}{1-2^{-2^{i+1}}}.\]
Hence, with Lemma \ref{varianteproduit} we get
\[\dfrac{\vert u_i\vert_0}{\vert u_i\vert}+\dfrac{\vert v_i\vert_1}{\vert v_i\vert}\geq2\times\dfrac{3}{4}\times\dfrac{1}{1-2^{-2^{i+1}}}\geq  \dfrac{3}{2}.
 \]
\end{itemize}
\end{proof}

\begin{lm}
 The letters of the  word $u$ do not admit unform frequencies.
\end{lm}

 \begin{proof} If the letters of $u$ possessed uniform frequencies, then the frequencies of 0 and 1, respectively denoted 
 $\mathbf{f}_u(0)$ and $\mathbf{f}_u(1)$, should verify  $\mathbf{f}_u(0)=\lim_{i\rightarrow \infty} \dfrac{\vert u_i\vert_0}{\vert u_i\vert}$, $\mathbf{f}_u(1)=\lim_{i\rightarrow \infty} \dfrac{\vert v_i\vert_1}{\vert v_i\vert}$ and $\mathbf{f}_u(0)+\mathbf{f}_u(1)=1$.  That is contradictory with Lemma \ref{excedent}.\end{proof}

\section{Complexity of $u$}

To estimate the complexity of $u$ we are going to observe its bispecial factors.

\begin{notation} \label{notation}
  Let $h,\, i\in \mathbb{N}$.  
We denote $u_{i}^{\left( h\right) }$ the finite word  $\sigma_h\sigma_{h+1}\sigma_{h+2}\ldots \sigma_{h+i-1}(0)$ and  $u^{\left( h\right) }$ the infinite word $\lim_{i\rightarrow \infty}u_{i}^{\left( h\right) }$. However  $u_{i}^{\left( 0\right) }$ and $u^{\left( 0\right) }$ are simply denoted respectively $u_i$ and $u$.
\end{notation}

\begin{definition}
 A  factor of $u^{\left( h\right) }$ is said to be short if it does not contain  $10$ as a factor. A  factor of $u^{\left( h\right) }$ which is not short is said to  be long.
\end{definition}
\begin{lm} (Synchronization lemma): Let $w$ be a long factor of $u^{\left( h\right) }$. Then there exist $x,\, y\in \mathcal{A}$ and $v\in \mathcal{A}^\ast$ such that $xvy$ is a factor of $u_{i}^{\left( h+1\right) }$ and $w=s\sigma_h\left( v\right) p$,  where $s$ is a non-empty suffix of $\sigma_h(x)$, and $p$ is a non-empty prefix of $\sigma_h(y)$. Moreover, the triple $(s,\,v,\,p)$ is unique. 
\end{lm}

\begin{proof} Since $w$ is long, it cannot occur inside the image of one letter. Any occurrence of $w$ in $u$ is therefore of the form $s\sigma_h(v)p$, so existence follows.

Uniqueness is consequence of the fact that $10$ occurs in $u^{\left( h\right) }$ only at the border between to images of letters under $\sigma_h$. \end{proof}

\begin{lm}
 \begin{enumerate}
  \item The short and strong bispecial factors  of $u^{(h)}$ are $\varepsilon$ and $1^{l_h}$.
  \item The short and weak bispecial factors  of $u^{(h)}$ are $0^{m_h-1}$ and $ 1^{n_h-1}$.
 \end{enumerate}

\end{lm}

 \begin{proof} Let us first observe that in $u^{\left( h\right) }$, the factor $01$ is always preceded by $10^{m_{h}-1}$. Therefore a bispecial factor containing $01$ must also contain $10$ and is long.

Then the short bispecial factors are all of the form  $0^k$ or $1^k$, $k\geq0$. We see that $\varepsilon$ is strong bispecial (extensions $00,\, 01,\, 10,\, 11$); $0^k$ ($0\leq k<m_h-1$) is ordinary bispecial (extensions $00^k0,\, 00^k1,\, 10^k0$), as well as $1^k$ ($1\leq k<n_h-1$, $k\neq l_h$); $1^{l_h}$ is strong bispecial (extensions $01^{l_h}0,\, 01^{l_h}1,\, 11^{l_h}0,\, 11^{l_h}1$); $0^{m_h-1}$ is weak bispecial (extensions $00^{m_h-1}1$ and $10^{m_h-1}0$), as well as $1^{n_h-1}$; $0^{m_h}$ and $1^{n_h}$ are not special, and $0^k$ ($k>m_h$) and $1^k$ ($k>n_h$) are not factors.
\end{proof}

\begin{lm}
 Let $w$ be a factor of  $u^{(h)}$. Then the following assertions are equivalent: \\
(1) $w$ is a long  bispecial factor of $u^{(h)}$.\\
(2) There exists a bispecial factor  $v$   of $u^{(h+1)}$ such that $w=\widehat{\sigma}_h(v)$ where $\widehat{\sigma}_h(v)=1^{l_h}\sigma_h(v)0^{m_h}1^{l_h}$.\\
Moreover $v$ and $w$ have the same type 
 and $|v|<|w|$.
\end{lm}

 \begin{proof}
First, let us observe this fact: If a finite word  $v$ is a factor of $u^{(h+1)}$ then $\widehat{\sigma}_h(v)=1^{l_h}\sigma_h(v)0^{m_h}1^{l_h}$ is a factor of $u^{(h)}$. Now, let us consider a bispecial factor $v$ of  $u^{(h+1)}$. Therefore the words $\widehat{\sigma}_h(0v),\ \widehat{\sigma}_h(1v),\ \widehat{\sigma}_h(v0)$ and $\widehat{\sigma}_h(v1)$ are factors of $u^{(h)}$; 
 moreover $0\widehat{\sigma}_h(v)$ and $1\widehat{\sigma}_h(v)$ are respectively  suffix of the first two words whereas  $\widehat{\sigma}_h(v)0$ and $\widehat{\sigma}_h(v)1$ are respectively prefix of the last two words. Hence, the word $w=\widehat{\sigma}_h(v)$ is bispecial in  $u^{(h)}$, and $\mathbf{m}(w)\geq \mathbf{m}(v).$

Conversely, let   $w$ be a long bispecial factor of $u^{(h)}$. Then, according to the synchronization lemma, we can write $w$  uniquely  in the form $s\sigma_{h}(v)p$ where $s$ and $p$ are respectively non-empty suffix and prefix of  images of letters. 

As $0w$ and $1w$ are factors of $u^{\left( h\right) }$, and $\sigma_h(v)p$ starts with $0$, it follows that $0s0$ and $1s0$ are factors of $u^{\left( h\right) }$. This is only possible if $s=1^{l_h}$ ( $s=1^k$ with $1\leq k<l_h$ or $l_h<k<n_h$ are excluded since $0s0\notin L\left( u^{\left( h\right) }\right) $;  $s=0^{k}1^{l_{h}}$ with $1\leq k<m_h$ and $s=0^{k}1^{n_{h}}$ with $0\leq k<m_h$ are excluded since $1s0\notin \mathcal{L}\left( u^{\left( h\right) }\right) $; and $s=0^{m_h}1^{l_{h}}$ and $s=0^{m_h}1^{n_{h}}$ are excluded since $0s0\notin \mathcal{L}\left( u^{\left( h\right) }\right) $).

Similarly, $1p0$ and $1p1$ are factors of $u^{\left( h\right) }$, and this is only possible if $p=0^{m_h}1^{l_{h}}$.
 Therefore $w=\widehat{\sigma}_h(v)$.

If $w$ extends as $awb$ with $a,\, b\in \mathcal{A}$, then $v$ also extends as $avb$. Therefore $\mathbf{m}(v)\geq \mathbf{m}(w).$ It follows that $\mathbf{m}(v)=\mathbf{m}(w)$: $v$ and $w$ have the same type. Moreover, it is clear that $|v|<|w|$.  \end{proof}

In fact, long bispecial factors of  $u^{(h)}$  are the  images by $\widehat{\sigma}_h$  of the ``less long'' bispecial factors of $u^{(h+1)}$. Thus,  step by step, any non-ordinary bispecial factor  $w$ of $u^{(h)}$ of given type, will be write in the following form $\widehat{\sigma}_h\widehat{\sigma}_{h+1}\ldots\widehat{\sigma}_{h+i-1}(v)$ where $v$ is a short bispecial factor of  $u^{(h+i)}$ with the same type.\\
We will call  bispecial factors of rank $i$, ($i\geq 0$) of $u^{(h)}$,  and write $a_{i}^{(h)}$, $b_{i}^{(h)}$, $c_{i}^{(h)}$,  $d_{i}^{(h)}$ the following words
 \[a_{i}^{(h)}=\widehat{\sigma}_h\widehat{\sigma}_{h+1}\ldots\widehat{\sigma}_{h+i-1}(\varepsilon),\ b_{i}^{(h)}=\widehat{\sigma}_h\widehat{\sigma}_{h+1}\ldots\widehat{\sigma}_{h+i-1}(1^{l_{h+i}}),\] \[c_{i}^{(h)}=\widehat{\sigma}_h\widehat{\sigma}_{h+1}\ldots\widehat{\sigma}_{h+i-1}(0^{m_{h+i}-1})\  \textrm{and}\  d_{i}^{(h)}=\widehat{\sigma}_h\widehat{\sigma}_{h+1}\ldots\widehat{\sigma}_{h+i-1}(1^{n_{h+i}-1}).\]
The short bispecial $\varepsilon,\ 1^{l_h}$, $0^{m_h-1}$ and $ 1^{n_h-1}$ of $u^{(h)}$ are the bispecial factors of rank $0$,  $a_{0}^{(h)}$, $b_{0}^{(h)}$, $c_{0}^{(h)}$, and  $d_{0}^{(h)}$.

The non-ordinary bispecial factors of $u$ are therefore $a_i=a_{i}^{(0)}$, $b_i=b_{i}^{(0)}$, $c_i=c_{i}^{(0)}$,  $d_i=d_{i}^{(0)}$.

\begin{definition} Let  $v, \, w\in \mathcal{A}^{\ast}$ and $V,\ W$ be their corresponding Parikh vectors. Let us  say that $V$ is less than  $W$ and write $V<W$ when $\vert v\vert_a\leq\vert w\vert_a$ for all $a\in \mathcal{A}$ and $\vert v\vert<\vert w\vert$. 

\end{definition}

\begin{prop} \label{compar1} Let   $v$, $w$, $v'$, $w'$ be four words such that $v'=\widehat{\sigma}_i(v)$ and $w'=\widehat{\sigma}_i(w)$. Then
\[ V<W\Longrightarrow V'<W' .\] 
\end{prop}

 \begin{proof} Assume that $V<W$. Then, $\vert v\vert_0\leq \vert w\vert_0$, $\vert v\vert_1\leq \vert w\vert_1$, and  $\vert v\vert<\vert w\vert$. On the one hand, we have $\vert v'\vert_0= m_i\left( \vert v\vert+1\right) $ and $\vert w'\vert_0=m_i\left( \vert w\vert+1\right) $; hence $\vert v'\vert_0<\vert w'\vert_0$.
On the other hand, we have $\vert v'\vert_1= l_i\vert v\vert_0+n_i\vert v\vert_1+2l_i $ and $\vert w'\vert_1= l_i\vert w\vert_0+n_i\vert w\vert_1+2l_i $; so $\vert v'\vert_1\leq\vert w'\vert_1$. Finally, $\vert v'\vert=\vert v'\vert_0+\vert v'\vert_1<\vert w'\vert_0+\vert w'\vert_1=\vert w'\vert_1.$  \end{proof}

\begin{lm} \label{vectordonne} For all $i\geq 0$, let $A_i,\, B_i,\, C_i,\,D_i$ be the Parikh vectors corresponding to the non-ordinary bispecial factors of $u$, $a_i,\, b_i,\, c_i,\,d_i$. 
 Then, we have
\[  
\forall i\geq 1,\ D_{i-1}< B_{i}< C_{i}< A_{i+1} < D_{i}
\]
\end{lm}

 \begin{proof}  
Applying $\widehat{\sigma}_{i-1} $ on the words  $b_{0}^{(i)}$, $c_{0}^{(i)}$, $\widehat{\sigma}_{i} \left( a_{0}^{(i+1)}\right) =1^{l_i}0^{m_i}1^{l_i}$, and  $d_{0}^{(i)}$ 
  we get the following words 
\[\left\{\begin{array}{lcl}
d_{0}^{(i-1)}&=& 1^{n_{i-1}-1}\\
b_{1}^{(i-1)}&=&1^{l_{i-1}}\left( 0^{m_{i-1}}1^{n_{i-1}}\right) ^{l_i}0^{m_{i-1}}1^{l_{i-1}}\\
c_{1}^{(i-1)}&=&1^{l_{i-1}}\left( 0^{m_{i-1}}1^{l_{i-1}}\right) ^{m_i-1}0^{m_{i-1}}1^{l_{i-1}}\\
a_{2}^{(i-1)}&=&1^{l_{i-1}}\left( 0^{m_{i-1}}1^{n_{i-1}}\right) ^{l_i}\left( 0^{m_{i-1}}1^{l_{i-1}}\right) ^{m_i}\left( 0^{m_{i-1}}1^{n_{i-1}}\right) ^{l_i}0^{m_{i-1}}1^{l_{i-1}}\\
d_{1}^{(i-1)}&=&1^{l_{i-1}}\left( 0^{m_{i-1}}1^{n_{i-1}}\right) ^{n_i-1}0^{m_{i-1}}1^{l_{i-1}}.\\
\end{array}\right. \]
The Parikh vectors corresponding to these words are: 
\[\left\{\begin{array}{lcl}
D_{0}^{(i-1)}&=&\begin{pmatrix}  0\\
                                  n_{i-1}-1 \\
                                 \end{pmatrix}\\
 
   B_{1}^{(i-1)}&=&\begin{pmatrix} m_{i-1}\left( l_i+1\right) \\
                                  n_{i-1}l_i+2l_{i-1}
                  \end{pmatrix}\\
   C_{1}^{(i-1)}&=&\begin{pmatrix} m_i m_{i-1}\\
                                  l_{i-1}\left( m_i+1\right) \\
                                 \end{pmatrix}\\
    A_{2}^{(i-1)} &=&\begin{pmatrix} m_{i-1}\left( m_i+2 l_i+1\right) \\
                                  l_{i-1}\left( m_i+2\right) +2l_in_{i-1}
                  \end{pmatrix} \\
D_{1}^{(i-1)}&=&\begin{pmatrix}  m_{i-1}n_i\\
                                  n_{i-1}\left( n_i-1\right)+2l_{i-1} \\
                                 \end{pmatrix}.\\
         \end{array} \right. \]

From  $ \left( \ast\right) $ we have $$n_{i-1}l_i+l_{i-1}<l_{i-1}m_i,\ m_i+2l_i+1<n_i, \ l_{i-1}m_i+2n_{i-1}l_i<n_{i-1}\left( n_i-1\right) .$$
It follows the inequalities:
\[  
  D_{0}^{(i-1)}<B_{1}^{(i-1)}< C_{1}^{(i-1)}< A_{2}^{(i-1)} < D_{1}^{(i-1)}
\]
Applying $\widehat{\sigma}_{i-2} $ on the words $d_{0}^{(i-1)}$, $b_{1}^{(i-1)}$, $c_{1}^{(i-1)}$, $a_{2}^{(i-1)}$, and  $d_{1}^{(i-1)}$ 
 we get the words $d_{1}^{(i-2)}$, $b_{2}^{(i-2)}$, $c_{2}^{(i-2)}$, $a_{3}^{(i-2)}$, and  $d_{2}^{(i-2)}$;
 By Proposition \ref{compar1},  it results the following inqualities:
\[  
 D_{1}^{(i-2)}< B_{2}^{(i-2)}< C_{2}^{(i-2)}< A_{3}^{(i-2)} < D_{2}^{(i-2)}.
\]
And so on, after the $i$-th iteration we get:
\[  
 D_{i-1}^{(0)}< B_{i}^{(0)}< C_{i}^{(0)}< A_{i+1}^{(0)} < D_{i}^{(0)}.
\]\end{proof}

\begin{lm} \label{bispordonne}
\[\forall i\geq 0, 
|b_{i}|<|c_{i}|< |a_{i+1}| < |d_{i}|< |b_{i+1}|
.\]
\end{lm}

 \begin{proof} $\bullet$ For $i\geq 1$,   the inequalities $
|b_{i}|<|c_{i}|< |a_{i+1}| < |d_{i}|< |b_{i+1}|
$ follows from Lemma \ref{vectordonne}\\
$\bullet$ For $i=0$, recall that  $$|b_0|=l_0,\ |c_0|=m_0-1,\  |a_1|=2l_0+m_0,\ |d_0|=n_0-1\ \textrm{and}\ |b_1|=l_1\left( m_0+n_0\right) +m_0+2l_0.$$ 
So
\[|b_{0}|<|c_{0}|< |a_{1}| < |d_{0}| < |b_{1}|
. \]\end{proof}

\begin{lm} \label{bispeciaux} The function $\mathbf{s} $ associated to the  word $u$ verifies:
 \[\forall n\in\mathbb{N},\ \mathbf{s}(n)=\left\{
\begin{array}{ll}
        1& \textrm{if}\  n=0\\
        2&  \textrm{if} \  n\in    \bigcup_{i\geq0}\Big( \left]|c_i|,\,|a_{i+1}| \right]\cup \left]|d_i|,\, |b_{i+1}| \right]\Big) \cup \left] 0,\,|b_0|\right]  \\
	3&  \textrm{if} \  n\in  \bigcup_{i\geq0}\Big( \left]|b_i|,\,|c_{i}| \right]\cup \left]|a_{i+1}|,\, |d_i| \right]\Big) .\\	
\end{array}\right.\]
\end{lm}

 \begin{proof} Let  $n\in\N$. We know that $a_i$, $b_i$, $c_i$, and $d_i$, $i\geq0$ are the only bispecial  factors of $u$ which are strong or weak. Hence, we have

$$\begin{array}{lr} \mathbf{s}\left(  n\right)  =& 1+\sum _{\begin{array}{c}w\ \textrm{bispecial}\\
                                                                             |w|<n
                                                             \end{array}}
\mathbf{m} \left( w\right)  \\
&=1+\#\left\lbrace  i\geq0:|a_{i}|<n\right\rbrace \\
&+\#\left\lbrace i\geq0:|b_{i}|<n\right\rbrace \\
&-\#\left\lbrace i\geq0:|c_{i}|<n\right\rbrace\\
&-\#\left\lbrace i\geq0:|d_{i}|<n\right\rbrace .

\end{array}$$

Since for $m\in]0,\, |b_0|[$ there is not strong or weak bispecial factor of $u$ with length $m$ 
we have, 
\[\textrm{for} \ 0<n\leq |b_0|,\ \mathbf{s}(n)=1+\sum _{\begin{array}{c}w\ \textrm{bispecial}\\
                                            |w|\leq n-1
                                           \end{array}}\mathbf{m} \left( w\right)=1+\mathbf{m}(\varepsilon)=2.\]
Suppose $n> |b_0|$. Then, there exists $i\in \N$ such that $n\in[|b_i|,\, |b_{i+1}|[$. Since the sequences $|a_i|$, $|b_i|$, $|c_i|$, and $|d_i|$ are increasing we are in one of the following cases:

$\bullet$  $n\in [|b_i|,\, |c_i|[$, then $\mathbf{s}(n)=1+ (i+1)+(i+1)-(i)-(i)=3.$

$\bullet$  $n\in [|c_i|,\, |a_{i+1}|[$, then $\mathbf{s}(n)=1+ (i+1)+(i+1)-(i+1)-(i)=2.$

$\bullet$  $n\in [|a_{i+1}|,\, |d_{i}|[$, then $\mathbf{s}(n)=1+ (i+2)+(i+1)-(i+1)-(i)=3.$

$\bullet$  $n\in [|d_i|,\, |b_{i+1}|[$, then $\mathbf{s}(n)=1+ (i+2)+(i+1)-(i+1)-(i+1)=2.$

\end{proof}

 \begin{thm} \label{complexite}
 
The complexity function  $\mathbf{p}$ of $u$ verifies:
 \[\forall n\geq1,\ \mathbf{p}(n)\leq3n+1.\]
\end{thm}

 \begin{proof} By  Lemma \ref{bispeciaux}, 
$\text{s}(n)\leq 3$ for all $n\geq 0$. So,
 $\mathbf{p}(n)=\mathbf{p}(0)+\sum^{n-1}_{m=0}\mathbf{s}(m)\leq \mathbf{p}(0)+3(n)=3n+1$. \end{proof}

\begin{prop} \label{compar2}
 Let $v$, $w$, $v'$, $w'$ be four finite words such that $v'=\widehat{\sigma}_i(v)$ and  $w'=\widehat{\sigma}_i(w)$. Then for all $\lambda>0$ we have:
\[ W>\lambda \left[ V+\begin{pmatrix}
                       1\\1
                      \end{pmatrix}\right] \Longrightarrow
 W'>\lambda \left[ V'+\begin{pmatrix}
                       1\\1
                      \end{pmatrix}\right]\]
\end{prop}

 \begin{proof} Assume that $W>\lambda \left[ V+\begin{pmatrix}
                       1\\1
                      \end{pmatrix}\right]$. Since $\vert v'\vert_0= m_i\left( \vert v\vert+1\right) $ and $\vert v'\vert_1= l_i\vert v\vert_0+n_i\vert v\vert_1+2l_i $ 
then: 
\[V'=\begin{pmatrix}
      m_i&m_i\\
     l_i&n_i
     \end{pmatrix}
\begin{pmatrix}
      |v|_0\\
      |v|_1
     \end{pmatrix} +\begin{pmatrix}
      m_i\\
     2l_i
     \end{pmatrix}.\] 
In the same way, we write $W'$ (it suffices to replace $V$ with $W$ in the previous formula).
It follows, 
\[W'-\lambda \left[ V'+\begin{pmatrix}
                       1\\1
                      \end{pmatrix}\right]=\begin{pmatrix}
      m_i&m_i\\
     l_i&n_i
     \end{pmatrix}
\begin{pmatrix}
      |w|_0-\lambda|v|_0\\
      |w|_1-\lambda|v|_1
     \end{pmatrix} +  \left( 1-\lambda\right) \begin{pmatrix}
      m_i\\
     2l_i
     \end{pmatrix} -\lambda\begin{pmatrix}
                       1\\1
                      \end{pmatrix}.\]
Since \[ W>\lambda V +\lambda\begin{pmatrix}
                       1\\1
                      \end{pmatrix}\ \textrm{and}\  \left( 1-\lambda\right) \begin{pmatrix}
      m_i\\
     2l_i
     \end{pmatrix}>-\lambda \begin{pmatrix}
      m_i\\
     2l_i
     \end{pmatrix} \] 
it follows that:
\[ W'-\lambda \left[ V'+\begin{pmatrix}
                       1\\1
                      \end{pmatrix}\right]>\lambda\left[ \begin{pmatrix}
      m_i&m_i\\
     l_i&n_i
     \end{pmatrix}\begin{pmatrix}
      1\\
      1
     \end{pmatrix}   - \begin{pmatrix}
      m_i+1\\
     2l_i+1
     \end{pmatrix} \right]=\lambda\begin{pmatrix}
                       m_i-1\\n_i-l_i-1
                      \end{pmatrix} >\begin{pmatrix}
                       0\\0
                      \end{pmatrix}.
 \]\end{proof}
This proposition allows to prove the following lemma:
\begin{lm} \label{large}
\[ \forall i\geq0,\ B_{i+1}>l_{i+1}\left[  D_{i}+\begin{pmatrix}
                       1\\1
                      \end{pmatrix}\right]
\]
\end{lm}

 \begin{proof} Let  us choose an integer $i\geq 1$. Then, we have  $b_{1}^{(i)}=\widehat{\sigma}_{i} \left( b_{0}^{(i+1)}\right) =1^{l_{i}}\left( 0^{m_{i}}1^{n_{i}}\right) ^{l_{i+1}}0^{m_{i}}1^{l_{i}}$ and  $d_{0}^{(i)}=1^{n_{i}-1}$; the  corresponding Parikh vectors are:
 $B_{1}^{(i)}=\begin{pmatrix}
      l_{i+1}m_i+m_i\\
     l_{i+1}n_i+2l_i
     \end{pmatrix}$ and $D_{0}^{(i)}=\begin{pmatrix}
      0\\
    n_i-1
     \end{pmatrix}$.
It follows the inequality:
\[ \ B_{1}^{(i)}>l_{i+1}\left[  D_{0}^{(i)}+\begin{pmatrix}
                       1\\1
                      \end{pmatrix}\right].
\]
 By regressive induction on $j\leq i$, suppose that:  
\[B_{i+1-j}^{(j)}>l_{i+1}\left[  D_{i-j}^{(j)}+\begin{pmatrix}
                       1\\1
                      \end{pmatrix}\right]\] 
where   $B_{i+1-j}^{(j)}$ and $ D_{i-j}^{(j)}$ are  respectively Parikh  vectors of the words $b_{i+1-j}^{(j)}$ and $ d_{i-j}^{(j)}$.\\
Thus, by Proposition \ref{compar2},
\[B_{i+2-j}^{(j-1)}>l_{i+1}\left[  D_{i-j+1}^{(j-1)}+\begin{pmatrix}
                       1\\1
                      \end{pmatrix}\right]\]
 since $B_{i+2-j}^{(j-1)}$ and  $ D_{i-j+1}^{(j-1)}$ are  respectively Parikh  vectors of $b_{i+2-j}^{(j-1)}=\widehat{\sigma}_{j-1} \left( b_{i+1-j}^{(j)}\right)$ and  $ d_{i-j+1}^{(j-1)}=\widehat{\sigma}_{j-1} \left( d_{i-j}^{(j)}\right).$
So, \[B_{i+1-j}^{(j)}>l_{i+1}\left[  D_{i-j}^{(j)}+\begin{pmatrix}
                       1\\1
                      \end{pmatrix}\right],\ 0\leq j\leq i.\]
 In the inequality above, we find the lemma by making $j=0$.\end{proof}

\begin{thm} The complexity function $\mathbf{p}$ of $u$ verifies 
 $\liminf \dfrac{\mathbf{p}(n)}{n}=2$
\end{thm}

 \begin{proof}  We have $\mathbf{s}(n)=2$  for $ |d_i|<n\leq|b_{i+1}|$. So 
\[\mathbf{p}\left( |b_{i+1}|\right) =\mathbf{p}\left( |d_{i}|\right) +2\left( |b_{i+1}|-|d_{i}|\right).\]
By  Lemma \ref{bispeciaux}, we have $\text{p}(n)\leq3n+1$ and we deduce that:
\[\mathbf{p}\left( |b_{i+1}|\right) \leq2 |b_{i+1}|+1 +\dfrac{1}{l_{i+1}}|b_{i+1}| \]
 since 
$B_{i+1}>l_{i+1}\left[  D_{i}+\begin{pmatrix}
                       1\\1
                      \end{pmatrix}\right]>l_{i+1}  D_{i}$.
So
 \[\dfrac{\mathbf{p}\left( |b_{i+1}|\right)}{ |b_{i+1}|}\leq 2+\dfrac{1}{|b_{i+1}|}+\dfrac{1}{l_{i+1}}\  \textrm{and} \ 
\lim_{i \rightarrow \infty}\dfrac{\mathbf{p}\left( |b_{i+1}|\right)}{ |b_{i+1}|}=2.\]
Thus, $\liminf \dfrac{\mathbf{p}(n)}{n}=2$, since $\mathbf{s}(n)\geq2$ (for all $n\geq 1$) implies 
 $\liminf \dfrac{\mathbf{p}(n)}{n}\geq2$.
 \end{proof}

\section{Proof of theorem \ref{rec.unif}}

Now, with Notation \ref{notation} we are able to explain the proof of Theorem \ref{rec.unif}.

 \begin{proof}  Let us show that for $i\geq 0$, there exists $N_i$ such that any factor of $u$ of length $N_i$ contains the prefix $u_i$. Indeed, $u$ does not contain $1^{n_0+1}$.

$\bullet$ For  $i=0$ any factor of $u$ of length $N_0=n_0+1$ contains the prefix $0=u_0$.

$\bullet$ For $i\geq1$,  any factor of $u^{(i)}$  of length $N_{0}^{\left( i\right) }=n_i+1$ contains the prefix $0=u_{0}^{(i)}$ of $u^{(i)}$.
Thus, any factor of $u^{(i-1)}$ of length   
$$N_{1}^{(i-1)}=\left( m_{i-1}+n_{i-1}\right) \left( N_0^{\left( i\right)} +1\right)  $$
contains $\sigma_{i-1}\left( 0\right)=u_1^{(i-1)}$.

By regressive induction on $j$, suppose that for $j\leq i-1$, there exists $N_{i-j}^{(j)}$ such that any factor of $u^{(j)}$ of length $N_{i-j}^{(j)}$ contains the word $u_{i-j}^{(j)}$. Then, any factor of $u^{(j-1)}$ of length
$$N_{i-j+1}^{(j-1)}=\left( m_{j-1}+n_{j-1}\right) \left( N_{i-j}^{(j)} +1\right)  $$
contains $\sigma_{j-1}\left(u_{i-j}^{(j)}\right) =u_{i-j+1}^{(j-1)}$.

 So, for $0\leq j\leq i-1$, there exists $N_{i-j}^{(j)}$  such that any factor  of $u^{(j)}$ of length $N_{i-j}^{(j)}$ contains the word $u_{i-j}^{(j)}$.

Consequently, letting $N_i=N_{i}^{(0)}$, it follows that any factor of $u=u^{(0)}$ of length  $N_{i}$ contains the word $u_i$. This completes the proof.\end{proof}

\textbf{Acknowledgmentts}

The authors would like to thank CNRS DERCI DSCA for its support during
this work.

\bigskip
 \begin{tabular}{l}
Julien CASSAIGNE\\
Institut de Math\'ematiques de Marseille\\
163 avenue de Luminy, case 907\\
F-13288 Marseille Cedex 9\\
France\\
julien.cassaigne@math.cnrs.fr
\end{tabular}

\bigskip
 \begin{tabular}{l}
Idrissa Kabor\'e\\
UFR-Sciences Exactes et Appliqu\'ees\\
Universit\'e Nazi Boni\\
01 BP 1091 Bobo-Dioulasso 01\\
Burkina Faso\\
ikaborei@yahoo.fr
\end{tabular}

\begin{thebibliography}{99}
\addcontentsline{toc}{chapter}{\protect\numberline{}{Bibliographie}} 

 \bibitem{bosher} M. Boshernitzan, \textit{A unique ergodicity of minimal symbolic flows with minimal block growth}, J. Analyse Math. 44 (1985), 77-96.
\bibitem{berthe} V. Berthé, \emph{Fréquences des facteurs des suites sturmiennes}, Theoret. Comput. Sci. 165 (1996), 295-309.
 \bibitem{cant} CANT, \emph{Combinatorics, Automata and Number Theory}, V. Berth\'e, M. Rigo (Eds),
Encyclopedia of Mathematics and its Applications 135,  Cambridge University Press (2010). 
\bibitem{cassaigne1} J. Cassaigne, \emph{Sequences with grouped factors}, in Developments in Language Theory III (DLT'97), pp. 211-222, Aristotle University of Thessaloniki, 1998.
\bibitem{cassaigne2} J. Cassaigne, \emph{Complexité et facteurs spéciaux}, Bull. Belg. Math. Soc.  4 (1997), 67--88.
\bibitem{cyrkra1} V. Cyr, B. Kra, \emph{Counting generic measures for a subshift of linear growth}, J. Eur. Math. Soc. (JEMS)  21 (2019), 355--380.
\bibitem{cyrkra2} V. Cyr, B. Kra, \emph{Realizing ergodic properties in zero entropy subshifts}, Isr. J.. Math.  240(1) (2020), 119--148.
\bibitem{dafi} M. Damron, J. Fickenscher,  \emph{The number of ergodic measures for transitive subshifts under the regular bispecial condition},Ergodic Theory Dyn. Syst. 42((1),  pp 86--140, 2022.
\bibitem{fogg} N. Pytheas Fogg, \emph{Substitutions in Dynamics, Arithmetics and Combinatorics}, Lectures Notes in Mathematics 1794, Springer-Verlag Berlin Heidelberg,  2002.
\bibitem{keane} M. Keane, \emph{Non-ergodic interval exchange transformations}, Israel J. Math. 26 (1977),  188--196.
\bibitem{monteil} T. Monteil, \emph{Illumination dans les billards polygonaux et dynamique symbolique}. Ph.D. thesis, Institut de Math\'ematiques de Luminy, Universit\'e de la M\'editerran\'ee, 2005. Chapter 5. 

\end{thebibliography}
\end{document}